\documentclass[11pt]{amsart}
\usepackage{verbatim}
\usepackage[utf8]{inputenc}
\usepackage{amsmath,amssymb,amsthm}
\usepackage{tikz-cd}
\usepackage{bbm}
\usepackage{hyperref}
\usepackage{comment}
\usepackage{stmaryrd}
\usepackage[new]{old-arrows}
\usepackage{graphicx}
\usepackage{todonotes}

\usepackage [english]{babel}
\usepackage [autostyle, english = american]{csquotes}

\newcommand{\lbm}{\left[ \begin{matrix}}
\newcommand{\rem}{\end{matrix} \right]}

\newtheorem{theorem}{Theorem}
\newtheorem{lemma}{Lemma}

\theoremstyle{definition}

\begin{document}
\title[Thomae's function and ergodic measures]{Thomae's function \\ and the space of ergodic measures}
\author[A.\ Gorodetski]{Anton Gorodetski}

\address{Department of Mathematics, University of California, Irvine, CA~92697, USA}

\email{asgor@uci.edu}

\author[A.\ Luna]{Alexandro Luna}

\address{Department of Mathematics, University of California, Irvine, CA~92697, USA}

\email{lunaar1@uci.edu}

\thanks{A.\ G.\ was supported in part by NSF grant DMS--1855541}

\date{}
\maketitle

\begin{abstract}
We study the space of ergodic measures of the map $$f:\mathbb{T}^2\to \mathbb{T}^2, \ f(x, y)=(x, \ x+y)(\text{mod}\, 1),$$ and show that its structure is similar to the graph of Thomae's function.
\end{abstract}
\section*{Introduction}
Thomae's function $\frak{T}:(0,1) \to \mathbb{R}$ given by
$$
\frak{T}(x)=\left\{
  \begin{array}{ll}
    0 & \hbox{if $x$ is irrational;} \\
    \frac{1}{q}, & \hbox{if $x=\frac{p}{q}$ for some coprime $p, q\in \mathbb{N}$}
  \end{array}
\right.
$$
and is also known as ``the popcorn function'', ``the raindrop function'', ``the countable cloud function'', the Riemann function, or even as ``the Stars over Babylon'' (as was suggested by John Horton Conway). Some of these names are justified by the form of the graph of this function, see Fig\ref{f.1}:
\begin{figure}[h]
\includegraphics[scale=0.5]{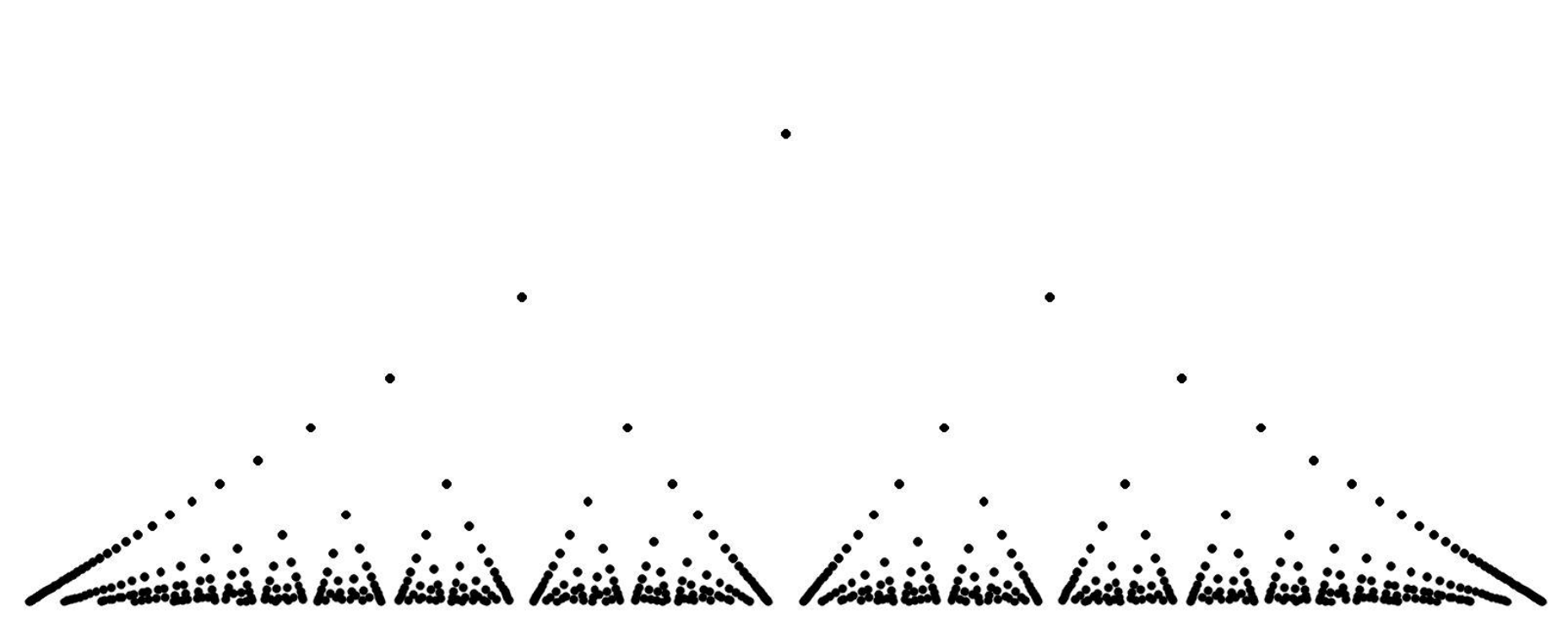}
\caption{The graph of Thomae's function.}\label{f.1}
\end{figure}

\

Thomae's function is a standard example that is presented in most introductory real analysis classes. This function is continuous at irrational values of the argument and is discontinuous at rational values. This concept is a starting point for an interesting journey towards the Baire Category Theorem, where one can prove that there does not exist a function which is continuous on the rationals and discontinuous on the irrationals. At the same time, Thomea's function itself remains a mathematical curiosity, and is something that doesn't usually appear in ``real life". The purpose of this short note is to show how a structure  similar to the graph of Thomae's function can appear in a natural way as a space of ergodic invariant measures of a very simple transformation.

\begin{center}
\begin{figure}[h]
\includegraphics[scale=.45]{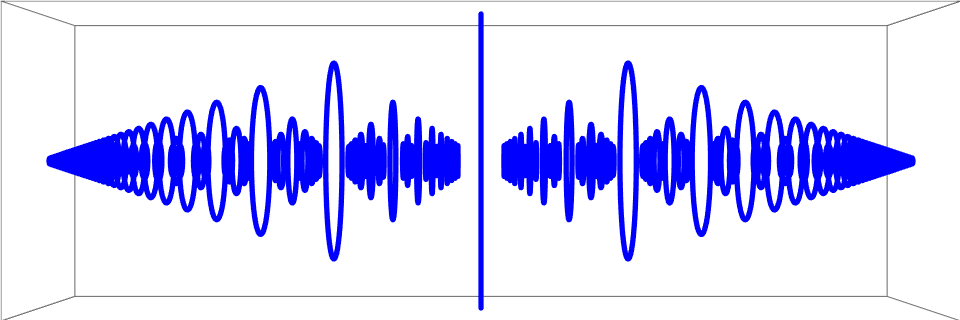}

\caption{The set $\mathcal{R}$ homeomorphic to the space of invariant ergodic measures of the map $f$.}
\label{f.2}
\end{figure}
\end{center}

Namely, define $f:\mathbb{T}^2\rightarrow \mathbb{T}^2$ via $f(x,\ y)=(x, \ x+y \ (\mathrm{mod} \ 1))$. Let $\delta_\omega$ be the atomic probability measure at the point $\omega\in \mathbb{T}^2$, and $\mathcal M(\mathbb{T}^2)$ be the space of all probability measures on $\mathbb{T}^2$. Consider the map  $T:\mathbb{T}^2\rightarrow \mathcal M \left(\mathbb{T}^2\right),$
$$T(\omega):=\lim\limits_{n\rightarrow\infty} \frac{\delta_{\omega}+\delta_{f(\omega)}+\dots+\delta_{f^{n-1}(\omega)}}{n},$$
where the limit is taken in the weak-$*$ topology. In this short note we describe the points at which $T$ is continuous, and show that its image has a topological structure similar to the graph of the Thomae's function. Namely, define the set $\mathcal{R}\subset \mathbb{R}^3$ that can be obtained as a revolution of the graph of the Thomae's function about $Ox$-axis, see Fig.\ref{f.2}:
$$
\mathcal{R}=\left\{(x,y,z)\in \mathbb{R}^3\ |\ x\in (0,1), \ \sqrt{y^2+z^2}=\mathfrak{T}(x)\right\}
$$
Then we have the following:
\begin{theorem}\label{Thm1}
The map $T$ is well defined, and has the following properties:
\begin{itemize}
\item   For any $\omega\in \mathbb{T}^2$ the image $T(\omega)$ is an ergodic invariant measure of the map $f$.
\item The map $T$ is continuous at Lebesgue almost every point of $\mathbb{T}^2$ and discontinuous at a dense set of points. In particular, $T$ is continuous at $\omega=(x_0, \ y_0)$ if and only if $x_0$ is irrational.
\item The image $T(\mathbb{T}^2)$ is homeomorphic to $\mathcal{R}$.
\end{itemize}
\end{theorem}
This result gives an example of a very simple transformation with a non-trivial topology of the space of ergodic measures. To put it into some context, let us mention that questions on the structure of spaces of invariant measures were studied in the theory of dynamical systems since 1970s.
K. Sigmund in \cite{S} 
showed that the space of ergodic measures for a transitive subshifts of finite type  is 
path-connected; the result was recently extended by I. Godofredo and V. Anibal \cite{GA} to the case of countable Markov shifts.  
Denseness of atomic ergodic measures in the space of all invariant measures for hyperbolic systems was studied by many authors, including 
 K. Sigmund \cite{S1} and 
 A. Katok \cite{KA}. 
 T. Downarowicz showed that for every Choquet simplex $K$ there exists a minimal subshift with the space of ergodic invariant measures affinely homeomorphic to $K$ \cite{D}.
 Questions  about structure of the space of hyperbolic ergodic measures on a locally maximal homoclinic class of a diffeomorphism were studied by A. Gorodetski and Y. Pesin \cite{GP}.
Some of the results by C. Bonatti and J. Zhang \cite{BZ}, by L.J. Diaz, K. Gelfert, T. Marcarini, and M. Rams  \cite{DGM}, and by D. Yang and J. Zhang  \cite{YZ} can be interpreted in terms of the structure of the space of ergodic invariant measures of partially hyperbolic diffeomorphisms under different assumptions. Recently, a notion of emergence have been developed by P. Berger and J. Bochi in \cite{BB}; it can be viewed as a way to quantify the ``size'' or ``dimension'' of the space of ergodic measures of a dynamical system. For other related results see also \cite{GK, KKK, KLW}.



\section*{Proof of the Main Result}

\noindent We first recall the following useful facts about sequences:
\begin{lemma}\label{sequence lemma}
\noindent
\begin{itemize}
    \item[(1)] If  $\left(\frac{p_n}{q_n}\right)_n$ is a sequence of fully reduced rationals with an irrational limit, then $q_n\rightarrow \infty$ as $n\rightarrow\infty$.
    \item[(2)] If $(a_n)_n$ is a sequence such that there is some $M$ satisfying $a_{n}=a_{n+M}$, then
    $$\lim\limits_{n\rightarrow\infty}\frac{a_1+\cdots+a_n}{n}=\frac{a_1+\cdots+a_M}{M}.$$
\end{itemize}
\end{lemma}
\noindent
Consider the ergodic invariant measures of $\mathbb{T}^2$ given by
$$\mu_{x_0}:=\delta_{x_0}\times m_{S^1} \ \ \text{and} \ \ \mu_{(p/q, \ y_0)}:=\frac{1}{q}\sum_{i=0}^{q-1}\delta_{\left(\frac{p}{q}, \ y_0+i\frac{p}{q}\right)}$$
where $m_{S_1}$ is the Lebesgue measure on $S^1$ and $p/q$ is a fully reduced rational.

\begin{lemma}\label{Range T}
For $\omega=(x_0, \ y_0)\in\mathbb{T}^2$,
$$T(\omega)= \begin{cases}
      \mu_{x_0} & \text{if} \ x_0\not\in\mathbb{Q} \\
       \mu_{(p/q, \ y_0)} & \text{if} \ x_0\in \mathbb{Q}, \ \text{with} \ x_0=\frac{p}{q} \ \text{fully reduced}
   \end{cases}.$$
\end{lemma}

\begin{proof} Let $\phi\in C_0\left(\mathbb{T}^2\right)$. For simplicity, denote
$$\mu_n:=\frac{\delta_{\omega}+\delta_{f(\omega)}+\dots+\delta_{f^{n-1}(\omega)}}{n}.$$
First suppose that $x_0\not\in\mathbb{Q}$. Then, we note that $y_0$ has dense orbit in $S^1$ under the rotation function $R_{x_0}$. Now, by the continuity of $\phi$ and the Weyl Equidistribution Theorem, we have
\begin{align*}
    \lim\limits_{n\rightarrow\infty} \int \phi \ d\mu_n &=\lim\limits_{n\rightarrow\infty} \sum_{i=0}^{n-1}\phi(x_0, R^i_{x_0}(y_0))\\
    &=\int_0^1 \phi \ d\left(\delta_{x_0}\times m_{S^1}\right).
\end{align*}
Now, suppose that $x_0=p/q$ is a fully reduced rational number.  Since the sequence $\left(\phi\left(\frac{p}{q}, \ y_0+i\frac{p}{q}\right)\right)_{i}$ satisfies the premises of Lemma \ref{sequence lemma} part (2), we have that $\mu_n$ converges to $\mu_{(p/q, \ y_0)}$ since
$$\lim\limits_{n\rightarrow \infty}\frac{1}{n} \sum_{i=1}^{n-1}\phi\left(\frac{p}{q}, \ y_0+i\frac{p}{q}\right)=\frac{1}{q}\sum_{i=0}^{q-1}\phi\left(\frac{p}{q}, \ y_0+i\frac{p}{q}\right).$$
\end{proof}
\begin{proof}[Proof of Theorem \ref{Thm1}]
That $T$ is well-defined and the first statement follow easily from the previous Lemma and construction of $T$.
We now prove the sufficient and necessary conditions for the continuity of $T$ at $\omega\in \mathbb{T}^2$.\\
\\
For sufficiency, it is enough to prove that if $\omega=\left(p/q, \ y_0\right)$, where $p/q$ is fully reduced, then $T$ is not continuous at $\omega$. We consider a continuous function $h:S^1\rightarrow S^1$ such that
$$\int_0^1h(y) \ dy>0 \ \text{and} \ h\left(y_0+i\frac{p}{k}\right)=0,$$
for $i=0, \ 1,\dots, \ q-1$. Now, let $((x_k, \ y_k))_k$ be a sequence such that $x_k\not\in\mathbb{Q}$ and $(x_k,\ y_k)\rightarrow (x_0, \ y_0)$, and consider the function $\phi\in C_0\left(\mathbb{T}^2\right)$ via $\phi(x, \ y)=h(y).$ Then, by Lemma \ref{Range T}, it is clear that
$$\lim\limits_{k\rightarrow\infty}\int_0^1 \phi \ d(T((x_k, \ y_k))\neq \int_0^1\phi \ d(T(\omega)).$$
Conversely, let $\omega=(x_0, \ y_0)\in\mathbb{T}^2$ where $x_0\not\in\mathbb{Q}$. We will show that if $\omega_k\rightarrow \omega$, then $T(\omega_k)\rightarrow T(\omega)$ in the weak-$^*$ topology. It suffices to consider the cases when $\omega_k=(x_k, \ y_k)$ with $x_k\not\in\mathbb{Q}$ and $\omega_k=(p_k/q_k, \ y_k)$, where each $p_k/q_k$ is a fully reduced rational value.\\
\\
Let $\phi\in C_0\left(\mathbb{T}^2\right)$.
The former case is almost immediate by the uniform continuity of $\phi$ and the Weyl Equidistribution Theorem.
In the latter case, we have that
$$\int_0^1\phi \ d(T(\omega_k))=\frac{1}{q_k}\sum_{i=0}^{q_k-1}\phi\left(\frac{p_k}{q_k}, \ y_k+i\frac{p_k}{q_k}\right).$$
As $\phi$ is uniformly continuous, for sufficiently large $k$, we can make the quantity $\phi\left(\frac{p_k}{q_k}, \ y_k+i\frac{p_k}{q_k}\right)$ arbitrarily close to $\phi\left(x_0, \ y_0+i\frac{p_k}{q_k}\right)$, so it suffices to show that
$$\lim\limits_{k\rightarrow\infty}\frac{1}{q_k}\sum_{i=0}^{q_k-1}\phi\left(x_0, \ y_0+i\frac{p_k}{q_k}\right)= \int\phi \ d(T(\omega)).$$
Notice that the partition of $S^1$ given by $P_k=\left\{R_{\frac{p_k}{q_k}}^i(y_0)\right\}_{i=0}^{q_k}$ splits $S^1$ into $q_k$ many intervals of equal length. Thus, by Lemma \ref{sequence lemma} part (1), since $q_k\rightarrow \infty$, we must have
\begin{align*}
    \lim\limits_{k\rightarrow\infty}\frac{1}{q_k}\sum_{i=0}^{q_k-1}\phi\left(x_0, \ y_0+i\frac{p_k}{q_k}\right)&=\lim\limits_{k\rightarrow\infty}\sum_{i=0}^{q_k-1}\phi\left(x_0, \ R_{\frac{p_k}{q_k}}^i(y_0) \right) \cdot (\text{mesh}(P_k))\\
    &=\int_0^1\phi(x_0, \ y) \ dy.
\end{align*}
To show that $T\left(\mathbb{T}^2\right)$ is homeomorphic to $\mathcal R$, we first show that it is homeomorphic to $\tilde{\mathcal R}= \left\{(x, \ \xi)\in S^1\times \mathbb{C}: \ |\xi|=\mathfrak{T}(x)  \right\}$, and then construct a homeomorphism from $\tilde{\mathcal R}$ to $\mathcal R$.\\
\\
Consider the map $\psi_1: T\left(\mathbb{T}^2\right)\rightarrow \tilde{\mathcal R}$ which sends $\mu_{x}\mapsto (x, \ 0)$, when $x\not\in\mathbb{Q}$, and $\mu_{(p/q, \ y)}\mapsto \left(\frac{p}{q}, \ \frac{1}{q} e^{2\pi i qy} \right)$ when $p/q$ is full reduced. Clearly $\psi_1$ is a bijection, and both $\psi_1$ and $\psi_1^{-1}$ 
are continuous.\\
\\
Set $\mathfrak T (0)=1$. To see that $\mathcal{\tilde R}$ is homeomorphic to $\mathcal R$, it is enough to show that the graph $G_{[0,1)}$ of $\mathfrak T$ on $[0, 1)$ is homeomorphic to the graph $G_{(0,1)}$ of $\mathfrak T$ on $(0,1)$. Consider the bijection $\psi_2: G_{[0,1)}\rightarrow G_{(0,1)}$ which maps $(0,\ 0)\mapsto \left(\frac{1}{2}, \frac{1}{2}\right)$, $\left(\frac{1}{n}, \frac{1}{n}\right)\mapsto \left(\frac{1}{n+1},\frac{1}{n+1}\right)$ for $n\in\mathbb{N}_{\geq 2}$, and is the identity otherwise.
It is clear that both $\psi_2$ and $\psi_2^{-1}$ are continuous, and hence $\psi_2$ is a homeomorphism.
\end{proof}


\end{document}